%% file: general.tex
\documentstyle[12pt,leqno]{article}

\input macros

\title{General existence of minimal surfaces of genus zero
with catenoidal ends and prescribed flux}
\date{}
\author{
   \begin{tabular}[t]{c}
    Shin Kato,\\
    {\footnotesize Osaka City mUniv.}
   \end{tabular}
   \begin{tabular}[t]{c}
    Masaaki Umehara\\
    {\footnotesize Osaka Univ.}
   \end{tabular}
   and 
   \begin{tabular}[t]{c}
    Kotaro Yamada\\
    {\footnotesize Kumamoto Univ.}
   \end{tabular}
}
\begin{document}

\maketitle

\input gen0


\input gen1


\input gen2


\input gen3


\input gena


\input genb


\input genr
\vspace{24pt}

\footnotesize
{\sc Shin Kato:} Department of Mathematics, Faculty of Science,
        Osaka City University, Osaka 558, JAPAN

\vspace{6pt}

{\sc Masaaki Umehara:} Department of Mathematics, 
        Graduate School of Science,
        Osaka University, Toyonaka 560, JAPAN
\par\noindent
{\it E-mail}: umehara@math.wani.osaka-u.ac.jp
\par\noindent

\vspace{6pt}

{\sc Kotaro Yamada:} Department of Mathematics, Faculty of Science,
        Kumamoto University, Kumamoto 860, JAPAN
\par\noindent
{\it E-mail:} kotaro@gpo.kumamoto-u.ac.jp
\end{document}

%% file: gen0.tex
\section*{Introduction}

Let 
$x:\hatC\setminus\{q_1,\ldots,q_n\}\to \R^3$ 
be a complete conformal minimal immersion. 
For each end $q_j$ $(j=1,\ldots,n)$ of $x$, 
the {\it flux vector\/} is defined by
\begin{equation}\label{eq01}
    \varphi_j := \int_{\gamma_j} \vec n \, ds, 
\end{equation}
where $\gamma_j$ is a positively oriented curve surrounding $q_j$,
and $\vec n$ the conormal such that 
$(\gamma', \vec n)$ is positively oriented. 
It is well known that 
the flux vectors satisfy a \lq\lq balancing" condition 
so called the {\it flux formula} 
\begin{equation}\label{eq02}
    \sum_{j=1}^n\varphi_j = 0. 
\end{equation}
\par
The minimal immersion $x$ is called an {\it $n$-end catenoid\/} 
if each end $q_j$ is of catenoid type. 
The catenoid and the Jorge-Meeks surfaces \cite{jm} are typical ones.
Recently, 
new examples of $n$-end catenoids have been found by 
\cite{kar}, \cite{l}, \cite{xu}, \cite{ross1}, \cite{ross2},
\cite{kat} and \cite{uy}.
For any $n$-end catenoid $x$, 
each flux vector $\varphi_j$ is proportional 
to the limit normal vector $\nu(q_j)$ 
with respect to the end $q_j$, 
and the scalar $w(q_j):=\varphi_j/4\pi\nu(q_j)$ 
is called the {\it weight\/} of the end $q_j$. 
In this case, 
the flux formula can be rewritten as follows. 
\begin{equation}\label{eq03}
    \sum_{j=1}^n 4\pi\,w(q_j)\,\nu(q_j) = 0. 
\end{equation}
It should be remarked that $w(q_j)$ may take a negative value. 
 \par
We consider the inverse problem
of the flux formula proposed in \cite{kat} and \cite{kuy}
as follows:

 \medskip

\noindent
{\bf Problem.} \quad 
{\it
    For given unit vectors $v:=\{v_1,\ldots,v_n\}$ in $\R^3$, 
    and nonzero real numbers $a:=\{a^1,\ldots,a^n\}$ 
    satisfying $\sum_{j=1}^n a^j v_j=0$ 
    $($we call such a pair $(v, a)$ flux data$)$, 
    is there a $($non-branched$)$ $n$-end catenoid 
    $x:\hatC\setminus \{q_1,\ldots,q_n\}\to \R^3$ 
    such that $\nu(q_j)=v_j$ 
    and $a_j$ is the weight at the end $q_j$?
}

 \medskip

We remark that 
Kusner also proposed a similar question (see \cite{ross1}).
Rosenberg and Toubiana \cite{rt} found  solutions with branch points 
in the category that the Gauss map is of degree $1$. 
But if one wishes a non-branched
solution, the degree of its Gauss map must be $n-1$, which
is the case just treated in this paper. 
 \par
The problem is not exactly affirmative. 
By the classification of Lopez~\cite{l}, 
we can see that 
the answer for $n\le 3$ is \lq\lq Yes" 
except for the case when two of $\{v_j\}_{j=1}^n$ coincide. 
Moreover, 
for $n\ge 4$,  
some obstructions exist as closed conditions 
in the space of flux data 
as shown in our previous paper \cite{kuy}. 
In spite of these obstructions, 
the authors also showed in \cite{kuy} that
the inverse problem is true 
for almost all flux data $(v, a)$
when $n=4$. 
In this paper, we treat the case $n\ge 5$ 
and show the following theorem: 

\vskip .5cm
\noindent
{\bf Theorem.}\quad 
{\it 
For each integer $n\ge 3$, 
the problem is solved for almost all flux data. 
}
\vskip .5cm

In Section 1, 
we reduce the inverse problem 
to seeking a sampling point 
satisfying certain non-degeneracy conditions. 
Two lemmas in Appendix A are applied to complete the reduction. 
In Section 2, 
we shall give a proof of Theorem. 
However, 
required technical calculations 
are done in Section 3 and Appendix B. 

The above general existence theorem does not
apply for the case that all flux vectors lie in the same plane,
since such flux data are contained in a measure zero subset in the 
set of all flux data. We say that such minimal surfaces
are of {\it Type II}. In [KUY2], we show that our approach in this paper
can be modified even for such a specified case and prove the
general existence of $n$-end catenoids ($n\le 8$) of Type II.
Recently, Kusner-Schmitt [KS] explain 
the moduli space of minimal surfaces with embedded planar ends 
by using the term of spin structure of Riemann surfaces.
It should be remarked that our approach can also be interpreted
in terms of spin structure. (See Remark~\ref{rmk:spin}.)

The author are very grateful to 
Professors Yusuke Sakane, Ichiro Enoki  and Koji Cho
for valuable discussions and encouragement.

%% file: gen1.tex
\section{Reduction}

The flux vector $\varphi_j$\,\,($j=1,...,n$)
 given by (\ref{eq01}) in introduction
can be rewritten as follows;
\begin{equation}
     \varphi_j  :=  
     -\mbox{Im}\left( 
         \oint_{\gamma_j}(1-g^2)\,\omega,
         \oint_{\gamma_j}\sqrt{-1}\,(1+g^2)\,\omega, 
         \oint_{\gamma_j} 2g\omega \right),
\end{equation}
where $(g,\omega)$ is the Weierstrass data of the
minimal immersion $x:\hatC\setminus\{q_1,\ldots,q_n\}\to \R^3$ 
given by
\begin{eqnarray*}
    g       &:=&  \partial x^3/(\partial x^1-\sqrt{-1}\,\partial x^2), \\
    \omega  &:=&  \partial x^1-\sqrt{-1}\,\partial x^2.
\end{eqnarray*}
On the other hand,
the well known Weierstrass representation is written as
\[
    x = 
      \mbox{Re}\left(
           \int_{z_0}^z(1-g^2)\,\omega,
           \int_{z_0}^z\sqrt{-1}\,(1+g^2)\,\omega, 
           \int_{z_0}^z 2g\omega \right).
\]
In particular, the monodromy vector of the immersion around 
the end $q_j$ (resp. the flux vector of $q_j$) is
the real part (resp. the imaginary part)
of the residue of the holomorphic vector
\[
    \partial x = 
    \frac{1}{2}
         \left(
           (1-g^2)\,\omega, 
           \sqrt{-1}\,(1+g^2)\,\omega,
           2g\omega
         \right), 
\]
around the end $z=q_j$.
We have shown  in the previous paper [KUY1] that the
inverse problem
of the flux formula reduces  to finding
solutions of a system of algebraic equations:

 \medskip

\begin{Thm}\,\,{\rm (\cite{kuy})} \label{Thm:1} 
Let $(v,a)$ be a pair of unit vectors 
$v=\{v_1,\ldots,v_n\}$ 
in $\R^3$ 
and nonzero real numbers 
$a=\{a^1,\ldots,a^n\}$ 
satisfying the balancing condition{\rm :} 
\begin{equation}\label{eq:bal}
    \sum_{j=1}^n a^j v_j=0.
\end{equation} 
Then there is an evenly branched $n$-end catenoid 
$x:\hatC\setminus\{q_1,\ldots,q_n\}\to \R^3$ 
$(q_j\ne\infty)$ 
such that the induced metric is complete at the end $q_j$, 
$\nu(q_j)=v_j$ 
and $a^j$ is the weight at the end $q_j$ 
$(j=1,\ldots,n)$, 
if and only if there exist complex numbers 
$b^1,\ldots,b^n$ 
satisfying the following conditions{\rm :} 
\begin{eqncases}{0.6}
    b^j\displaystyle{\Condsum{k=1}{k\ne j}^nb^k\frac{p_j-p_k}{q_j-q_k}}
        = a^j & \label{eq:y} \\[-1.5ex]
    &\qquad\qquad(j=1,\dots,n),\nonumber\\[-1.5ex]
    b^j\displaystyle{\Condsum{k=1}{k\ne j}^n 
                     b^k\frac{\overline{p_j}p_k+1}{q_j-q_k}}
        = 0 &
          \label{eq:x}
\end{eqncases}
where $p_j:=\sigma(v_j)$, 
$\sigma:S^2\to\hatC$ is the stereographic projection, 
and we assume $p_j\ne\infty$. 
 \par
Moreover, 
the surface $x$ has no branch points 
if and only if the two the polynomials 
\begin{eqnarray}
    Q(z) & := \displaystyle\sum_{j=1}^nb^j
              \displaystyle\Condprod{k=1}{k\ne j}^n(z-q_k), 
              \label{Def:Q} \\
    P(z) & := \displaystyle\sum_{j=1}^np_jb^j
              \displaystyle\Condprod{k=1}{k\ne j}^n(z-q_k) 
              \label{Def:P} 
\end{eqnarray}
are mutually prime and one of them has degree $n-1$. 
\end{Thm}

 \medskip

\begin{Rmk} \label{rmk:1} 
When $p_j=rq_j$, 
the theorem reduces to the results in the first author \cite{kat}. 
In this case the system (\ref{eq:y}) and (\ref{eq:x}) reduces to 
$$
    \left\{
    \begin{array}{l}
        rb^j \displaystyle{\Condsum{k=1}{k\ne j}^n} b^k 
         = a^j \\ 
        b^j\displaystyle\Condsum{k=1}{k\ne j}^n 
         b^k \displaystyle{\frac{|r|^2\overline{q_j}q_k+1}{q_j-q_k}} 
         = 0 
    \end{array}
    \right.
    \qquad\qquad 
    (j=1,\ldots,n). 
$$
As seen in \cite{kat}, 
the surface has no branch point 
if and only if $\beta:=\sum_{j=1}^n b^j\ne 0$. 
By using the relation 
$
P(z)/Q(z)=rz-r\beta/\left(\sum_{j=1}^n b^j/(z-q_j)\right) 
$, 
it is also checked directly 
from the last condition of the theorem. 
\end{Rmk}

 \medskip

\begin{Rmk} \label{rmk:2} 
The position of the ends 
$\{q_1,\ldots,q_n\}$ in the source domain $\hatC$ 
has the freedom of M\"obius transformations. 
For example, 
the following normalization is possible: 
$$
 q_1=1,\qquad q_{n-1}+q_{n-2}=0,\qquad q_n=0. 
$$
\end{Rmk}

 \medskip

\begin{Rmk} \label{rmk:3} 
The system of the equations 
({\rm \ref{eq:y}}) and ({\rm \ref{eq:x}}) 
has another expression 
\begin{eqncases}{0.6}
    b^j\displaystyle{\Condsum{k=1}{k\ne j}^nb^k\frac{1}{q_j-q_k}} 
     & = \displaystyle{a^j\frac{\overline{p_j}}{|p_j|^2+1}}, 
     \label{eq:xx}\\ 
    b^j\displaystyle{\Condsum{k=1}{k\ne j}^nb^k\frac{p_j+p_k}{q_j-q_k}} 
     & = \displaystyle{a^j\frac{|p_j|^2-1}{|p_j|^2+1}}.
     \label{eq:yy}
\end{eqncases}
Moreover we may replace (\ref{eq:xx}) by
\begin{equation}
   p_j b^j\displaystyle{\Condsum{k=1}{k\ne j}^nb^k\frac{p_k}{q_j-q_k}} 
     = -\displaystyle{a^j\frac{p_j}{|p_j|^2+1}}.
     \label{eq:xxx}\\
\end{equation}
In fact, 
if we set 
$$
    \gamma_j := b^j\Condsum{k=1}{k\ne j}^nb^k\frac{1}{q_j-q_k}, 
    \quad
    \delta_j := b^j\Condsum{k=1}{k\ne j}^nb^k\frac{p_k}{q_j-q_k} 
    \qquad (j=1,\ldots,n), 
$$
then (\ref{eq:y}) and (\ref{eq:x}) are written as
$$
    p_j\gamma_j-\delta_j = a^j, 
    \qquad 
    \gamma_j+\overline{p_j} \delta_j = 0. 
$$
It is equivalent to the relations 
$$
    \gamma_j = a^j \frac{\overline{p_j}}{|p_j|^2+1}, 
    \qquad 
    p_j\gamma_j+\delta_j = a^j \frac{|p_j|^2-1}{|p_j|^2+1}, 
$$
that is (\ref{eq:xx}) and (\ref{eq:yy}). 
On the other hand, 
$$
    p_j\gamma_j 
     = a^j\frac{|p_j|^2}{|p_j|^2+1} 
     = a^j\frac{|p_j|^2-1}{|p_j|^2+1}+\frac{a^j}{|p_j|^2+1} 
     = p_j\gamma_j+\delta_j+\frac{a^j}{|p_j|^2+1}, 
$$
which yields (\ref{eq:xxx}). 
\end{Rmk}

 \medskip

\begin{Rmk}\label{rmk:spin}
  The construction of $n$-end catenoids mentioned above is related to
  the spinor representation of minimal surfaces (cf.~\cite{ks});
  \[
     x = \mbox{Re}\left(
              \int_{z_0}^z (s_1{}^2-s_2{}^2),
              \int_{z_0}^z \sqrt{-1}\,(s_1{}^2+s_2{}^2),
              \int_{z_0}^z 2s_1\,s_2
            \right),
  \]
  where $(s_1,s_2)$ is a pair of meromorphic sections of the
  half-canonical bundle on $\C\cup\{\infty\}$.
  In fact, $s_1$ and $s_2$ have the following explicit expressions
  in this case:
  \[
    s_1 := \frac{Q(z)}{R(z)}\sqrt{-dz},\qquad
    s_2 := \frac{P(z)}{R(z)}\sqrt{-dz},
  \]
  where we set 
  $
     R(z):=\displaystyle\prod_{k=1}^n (z-q_k).
  $
\end{Rmk}

 \medskip

Theorem \ref{Thm:1} produces many $n$-end catenoids 
as seen in \cite{kat} and \cite{kuy}. 
First, 
we fix our attention to the equation (\ref{eq:x}). 
We consider a matrix 
\begin{equation}
    A_p := \left(\frac{\overline{p_j}p_k+1}{q_j-q_k}
           \right)_{j,k=1,\ldots,n}, 
    \label{A}
\end{equation}
where the diagonal components are interpreted as $0$. 
Then the vector ${}^t(b^1,\ldots,b^n)$ 
belongs to the kernel of the matrix $A_p$. 
As shown in the later sections, 
it is reasonable to expect that 
the rank of the matrix $A_p$ is generically $n-1$. 
In this case, 
${}^t(b^1,\ldots,b^n)$ should be proportional to 
any column vector of the cofactor matrix $\widetilde A_p$ of $A_p$. 
(By the definition, 
$A_p\widetilde A_p = \widetilde A_p A_p =(\det A_p)I$ 
holds.) 
So we set 
$$
    b_p(q) = {}^t(b^1_p(q),\ldots,b^n_p(q)) 
          := \mbox{the $n$-th column of the cofactor matrix} 
             ~ \widetilde A_p(q). 
$$

Now we projectify the problem: 
For fixed $p:=(p_1,\ldots,p_n)\in \C$, 
define a rational map between two complex projective spaces
$$
    \Fl_p = [f^1_p,\ldots,f^n_p]:\Pj^{n-1}\to \Pj^{n-1} 
$$
by 
\begin{equation} \label{eq:F} 
    f^j_p(q_1,\ldots,q_n) 
     := b_p^j(q)\sum_{k\ne j}b_p^k(q)
                \frac{p_j-p_k}{q_j-q_k}
        \qquad (j=1,\ldots,n).
\end{equation}
We set
$$
    \fl^j_p(q):=\Delta(q)^5\cdot f_p^j(q), 
$$
where
 $\Delta(q)$ is the difference product 
defined by 
\begin{equation}
    \Delta(q_1,\ldots,q_{n}) := \prod_{j>k}^n(q_j-q_k). 
\end{equation}
 \par
It is easily seen that 
each $\fl^j_p$ is a homogeneous polynomial in $q_1,\ldots,q_n$
and  $\Fl_p$ has another expression
$$
    \Fl_p=[\fl^1_p,\ldots,\fl^n_p]. 
$$
This projective formulation is reasonable 
in the following two senses: 
\begin{itemize}
\item 
Any homothety of $n$-end catenoids changes 
their weights $(a^1,\ldots,a^n)$ only by a constant multiplication. 
It allows us to projectify the image of $\Fl_p$. 
\item 
Changing coordinates of $n$-end catenoids 
by homothetic transformations 
corresponds to complex multiplications of $(q_1,\ldots,q_n)$. 
(See Remark \ref{rmk:2}.) 
It allows us to projectify the domain of $\Fl_p$. 
\end{itemize}
\par
Since $p_j$ is the stereographic image of $v_j$, 
the balancing condition (\ref{eq:bal}) is rewritten as 
$$
     \sum_{j=1}^n \frac{|p_j|^2-1}{|p_j|^2+1} a^j = 0, \qquad
     \sum_{j=1}^n \frac{\overline{p_j}}{|p_j|^2+1} a^j = 0. 
$$
We define a subspace $\W_p^{n-4}$ in $\Pj^{n-1}$ by 

\footnotesize

\begin{eqnarray*}
    \W_p^{n-4} 
      :=    \left\{[a^1,\ldots,a^n]\in \Pj^{n-1} \,\,; 
            \sum_{j=1}^n \frac{|p_j|^2-1}{|p_j|^2+1} a^j = 0, 
            \sum_{j=1}^n \frac{\overline{p_j}}{|p_j|^2+1} a^j = 0, 
            \sum_{j=1}^n \frac{p_j}{|p_j|^2+1} a^j = 0 
            \right\}. 
\end{eqnarray*}

\normalsize

We will show that for open dense $p\in \C^n$, 
the image of the map $\Fl_p$
is open dense in $\W_p^{n-4}$, and next show that 
it covers open dense subset of the totally real set
$\W_{\R}=\{[a]\in \W_p^{n-4}\,;\, a_j\in \R\}$.
Then the image of the map $\Fl_p$ 
contains $[a]\in \W_\R$ for almost all flux data $(p,a)$,
and Theorem in Introduction is obtained.
If $\Fl_p$ is a holomorphic map and there is a 
point at which the rank of $d\Fl_p$ is
$n-4$, the surjectivity of the map follows by the proper mapping theorem.
(See \cite{gr}.)
But unfortunately, the map $\Fl_p$ is singular on the set 
$\bigcap_{j=1}^n \Z(\fl^j_p)$, 
where $\Z(\fl^j_p)$ is the set of zeros of $\fl^j_p$.
As shown below,
we will overcome this difficulty by a usual blowing up process. 
 
From here, 
assume $\dim\langle v_1,\ldots,v_n\rangle=3$, 
where $v_j:=\sigma^{-1}(p_j)$ 
and $\sigma$ is the stereographic projection. 
Then clearly $\dim\W_p^{n-4}=n-4$. 
We remark here that 
$\dim\W_p^{n-4}=n-4$ 
holds for open dense $p\in\C^n$. 
Now we have the following lemma:

 \medskip

\begin{Lemma} \label{lamdalem}
For each $p\in \C^n$, 
the following relation holds: 
$$
    \Fl_p\biggl(\Z(\lambda_p)\setminus
                \bigcap_{j=1}^n \Z(\fl^j_p)\biggr)
     \subset \W_p^{n-4}, 
$$
where $\lambda_p$ is the determinant 
of the matrix $\Delta\cdot A_p$ 
and $\Z(\lambda_p)$ is the set of zeros 
of the homogeneous polynomial $\lambda_p$. 
\end{Lemma}

\begin{Proof}
Let 
$q\in\Z(\lambda_p)\setminus\bigcap_{j=1}^n \Z(\fl^j_p)$. 
If $\Delta(q)=0$, 
then it is easy to see that 
$q\in\bigcap_{j=1}^n \Z(\fl^j_p)$. 
Hence $\Delta(q)\ne 0$, 
and we get (\ref{eq:y}) with $b^j=b^j(q)$ $(j=1,\ldots,n)$. 
Recall Remark \ref{rmk:3}. 
Then the assertion of the lemma immediately follows 
by summing up (\ref{eq:yy}), 
(\ref{eq:xx}) and (\ref{eq:xxx}) for $j=1,\ldots,n$. 
\end{Proof}

 \medskip

We define an $(n-1)$-matrix $J_p$ by 
\begin{equation}
    J_p := \left((f^n_p)^2 
            \left\{\frac{\partial\det A_p}{\partial q_n}
              \cdot 
             \frac{\partial \stackrel\circ{f_p^k}}{\partial q_j} 
            -\frac{\partial\det A_p}{\partial q_j} 
              \cdot 
	       \frac{\partial \stackrel\circ{f_p^k}}{\partial q_n}\right\}
           \right)_{k,j=1,\ldots,n-1}, 
\end{equation}
where 
$$
    \stackrel\circ{f^j_p}:=\frac{f^j_p}{f^n_p} \qquad (j=1,\ldots,n-1).
$$
The matrix $J_p$ has a direct expression

\footnotesize

$$
    \mbox{\normalsize $J_p$} = \left(
              \frac{\partial\det A_p}{\partial q_n}\cdot 
              \left\{\frac{\partial f_p^k}{\partial q_j}\cdot f_p^n 
              -f_p^k\cdot\frac{\partial f_p^n}{\partial q_j}\right\} 
             -\frac{\partial\det A_p}{\partial q_j}\cdot 
              \left\{\frac{\partial f_p^k}{\partial q_n}\cdot f^n_p 
              -f^k_p \cdot \frac{\partial f_p^n}{\partial q_n}\right\}
           \right)_{k,j=1,\ldots,n-1}. 
$$

\normalsize

The following proposition plays an important role 
to establish Theorem in Introduction. 

 \medskip

\begin{Prop} \label{Prop:A} 
Suppose that there exist $u_0\in \C^n$ and 
a point $c=[c_1,\ldots,c_n]\in \Pj^{n-1}$
satisfying the following conditions{\rm :}
\begin{enumerate}
\item[\rm (1)]
  $c_1,\ldots,c_n$ are all distinct{\rm;}
\item[{\rm (2)}] 
    The rank of the matrix $A_{u_0}(c)$ is $n-1${\rm ;} 
\item[\rm (3)]
$\frac{\partial\det A_{u_0}}{\partial q_n}$ 
does not vanish at $q=c${\rm ;}
\item[{\rm (4)}] 
    The rank of the matrix $J_{u_0}(c)$ is $n-4${\rm ;} 
\item[{\rm (5)}] 
    Two polynomials $P(z)$ and $Q(z)$ 
    defined in\/ {\rm (\ref{Def:P})} and\/ {\rm (\ref{Def:Q})} 
    associated with the data $(q,p)=(c,u_0)$ 
    and $b=b_{u_0}(c)$ 
    are mutually prime and 
    one of them has degree $n-1${\rm ;}
\item[{\rm (6)}] 
    $f_{u_0}^j(c)\ne 0$ $(j=1,\ldots,n)${\rm ;} 
\item[{\rm (7)}] 
    $c_j\ne 0$ $(j=1,\ldots,n-1)$. 
\end{enumerate}
Then there exists an open dense subset $U\subset \C^n$
and an open dense subset $\Omega_p$ of the totally real set 
$\W_{\R}=\{[a]\in \W_p^{n-4}\,;\, a_j\in \R\}$
such that, for any $p\in U$ and $[a]\in\Omega_p$,
there exists an {\rm (}non-branched{\rm )} $n$-end catenoid with the 
flux data $(p,a)$.
\end{Prop}

 \medskip

By the proposition, 
the inverse problem of the flux formula can be solved 
for almost all flux data 
if one succeeds to take such a point $c$. 
This will be done in the next section. 
The outline of the proof of the proposition is as follows. 

 \medskip

By the condition (4), 
at least one $(n-4)$-submatrix $S_{u_0}$ of $J_{u_0}$ 
is of rank $n-4$. 
Let $1\le j_1<j_2<\cdots < j_{n-4}< n$ be  
the indices of the columns of the submatrix $S_{u_0}$, 
and $\{m_1,m_2,m_3\}$ their complement, 
namely $\{m_1,m_2,m_3\}=\{1,\ldots,n-1\}\setminus\{j_1,\ldots,j_{n-4}\}$. 
By Remark \ref{rmk:2}, 
we may restrict the flux map into the
following subspace of $\Pj^{n-1}$ 
containing the sampling point $c$: 
$$
    \V^{n-3} := \{[q_1,\ldots,q_n]\in \Pj^{n-1} \, ; \, 
                c_{m_2}q_{m_1}-c_{m_1}q_{m_2}=0, \,\, 
                c_{m_3}q_{m_1}-c_{m_1}q_{m_3}=0\}. 
$$
Now we define a homogeneous polynomial in $q_1,\ldots,q_n$ by 
$$
    H_{p}(q) :=   \Delta(q)^2\frac{\partial\det A_p}{\partial q_n}(q) \cdot 
                  \det\left(\Delta(q)^\ell S_{p}(q)\right) \cdot 
                  R_{p}(q) \cdot 
                  \prod_{j=1}^n\fl_{p}^j(q) \cdot 
                  \prod_{k=1}^{n-1}q_k, 
$$
where $\ell$ is chosen sufficiently large so that
$\det(\Delta(q)^\ell S_{p}(q))$ is a homogeneous polynomial
in $q_1,\ldots,q_n$, and 
$R_{p}$ is the resultant of the two polynomials 
$P(z)$ and $Q(z)$ of degree $n-1$ 
defined by (\ref{Def:P}) and (\ref{Def:Q}). 
(It can be easily shown that 
$R_p$ is also a homogeneous polynomial 
with respect to $q$. 
Or one may replace $R_q$ 
by the resultant of $P(q_1z)$ and $Q(q_1z)$.) 
Then by the conditions (1)-(7), 
$c\in \V^{n-3}$ satisfies $H_{u_0}(c)\ne 0$. 
We prove the following 

 \medskip

\begin{Lemma} \label{lem:1} 
The subset
$$
    U:=\{p\in \C^n\,;\, \Z(\lambda_p)\cap \V^{n-3}\not\subset 
    \Z(H_p)\}
$$
is open dense in $\C^n$,
where $\lambda_p=\det(\Delta\cdot A_p)$ is the homogeneous polynomial
defined in Lemma \ref{lamdalem}.
\end{Lemma}

\begin{Proof}
Obviously 
$U$ is an open subset of $\C^n$. 
Suppose that $U$ is not dense in $\C^n$. 
Then there exists an open subset $V$ such that 
\begin{equation}
    \Z(\lambda_p|_{\V^{n-3}})\subset \Z(H_p|_{\V^{n-3}}) 
     \qquad (p\in V). 
     \label{eq:Z}
\end{equation}
Since $\V^{n-3}\cong \Pj^{n-3}$, 
by Lemma A.1 in Appendix, 
(\ref{eq:Z}) holds for any $p \in \C^n$ 
such that $\lambda_p\not\equiv 0$. 
But this contradicts the fact that 
$\lambda_{u_0}(c)=0$, 
$\lambda_{u_0}\not\equiv 0$ 
and $H_{u_0}(c)\ne 0$. 
\end{Proof}

 \medskip

Roughly speaking, 
if $\Fl_p$ has no singularities and is of maximal rank, 
then it is surjective and 
we can find a pair $(q,b_p(q))$ 
satisfying (\ref{eq:y}) and (\ref{eq:x}). 
But unfortunately, 
$\Fl_p$ has singularities 
on $\bigcap_{j=1}^n \Z(\fl^j_p)$.
For this reason, 
we define a new variety $\hat\V^{n-3}$ 
and a map $\widehat{\Fl_p}\colon{}\hat\V^{n-3}\to\W_p^{n-4}$ 
instead of $\V^{n-3}$ and $\Fl_p$ as follows.
First we consider an algebraic set
\begin{eqnarray*}
    \Y^{n-3} 
     & = & \Biggl\{([q_1,\ldots,q_n],[a^1,\ldots,a^n]) 
           \in \Pj^{n-1}\times \Pj^{n-1} \,\, ; \\ 
     &   & \qquad 
           c_{m_2}q_{m_1}-c_{m_1}q_{m_2}=0, \,\, 
           c_{m_3}q_{m_1}-c_{m_1}q_{m_3}=0, \\ 
     &   & \qquad 
           a^j\fl_p^k=a^k\fl_p^j \qquad (j,k=1,\ldots,n), \\ 
     &   & \qquad 
           \sum_{j=1}^n \frac{|p_j|^2-1}{|p_j|^2+1} a^j = 0, 
           \sum_{j=1}^n \frac{p_j}{|p_j|^2+1} a^j = 0, 
           \sum_{j=1}^n \frac{\overline{p_j}}{|p_j|^2+1} a^j = 0 \\ 
     &   & \hspace{6cm} (j=1,\ldots,n) \Biggr\}, 
\end{eqnarray*}
and define two canonical projections: 
\begin{eqnarray*}
    & & \pi :\Y^{n-3}\ni ([q],[a]) \mapsto [q] \in \V^{n-3}, \\
    & & \pi':\Y^{n-3}\ni ([q],[a]) \mapsto [a] \in \W_p^{n-4}. 
\end{eqnarray*}
These two projections are both well-defined on $\Y^{n-3}$. 
Let $\hat\V^{n-3}$ be the algebraic closure of the set 
\begin{equation}
    \hat\V_{\mbox{\scriptsize reg}}^{n-3} 
     := \pi^{-1}\left(\V^{n-3} \setminus 
        \bigcap_{j=1}^n \Z(\fl^j_p)\right). 
        \label{set}
\end{equation}
We denote 
the restriction of the first projection $\pi$ to $\hat \V^{n-3}$ 
also by $\pi$.
We remark that 
$\pi|_{\hat\V_{\mbox{\scriptsize reg}}^{n-3}}:
\hat\V_{\mbox{\scriptsize reg}}^{n-3}\to 
\V^{n-3}\setminus\bigcap_{j=1}^n \Z(\fl^j_p)$ 
is bijective. 
On the other hand, 
we denote 
the restriction of the second projection $\pi'$ to $\hat \V^{n-3}$ 
by 
$$
    \widehat{\Fl_p}:\hat \V^{n-3}\to \W^{n-4}_p. 
$$
The map $\Fl_p\circ\pi$ is well-defined 
on $\hat\V_{\mbox{\scriptsize reg}}^{n-3}$, 
and coincides with the map $\widehat{\Fl_p}$.

 \medskip

\begin{Lemma} \label{lem:2} 
For each $p\in U$ 
satisfying $\dim\W_p^{n-4}=n-4$, 
there exists an irreducible component $\hat X^{n-4}$ 
of the algebraic set $\Z(\lambda_p\circ \pi)\cap \hat\V^{n-3}$ 
such that $H_p\circ \pi$ 
is not identically zero on $\hat X^{n-4}$. 
In addition, 
the restriction of the lifted flux map 
$\widehat{\Fl_{p}}|_{\hat X^{n-4}}:\hat X^{n-4}\to \W^{n-4}_p$ 
is surjective. 
\end{Lemma}

\begin{Proof}
Suppose that 
$\Z(\lambda_p\circ \pi)\cap \hat\V^{n-3}\subset \Z(H_p\circ \pi)$. 
Since $H_p$ is identically zero on the singular set 
$\bigcap_{j=1}^n \Z(\fl^j_p)$, 
it follows that 
$$
    \Z(\lambda_p)\cap \V^{n-3}\subset \Z(H_p).
$$
But this contradicts Lemma \ref{lem:1}. 
Hence 
there exists an irreducible component $\hat X^{n-4}$ 
of the algebraic set $\Z(\lambda_p\circ \pi)\cap \hat\V^{n-3}$ 
such that $H_p\circ \pi$ 
is not identically zero on $\hat X^{n-4}$. 
We set 
$$
    X^{n-4} := \pi(\hat X^{n-4}).
$$
Now we take a point $x_0\in X^{n-4}$ such that $H_p(x_0)\ne 0$. 
Consequently, 
we have $x_0\not\in \bigcap_{j=1}^n \Z(\fl^j_p)$ 
and so $\Fl_p(x_0)\in \W_p^{n-4}$ exists. 
We remark here that 
$m_1$-th component of $x_0$ in the homogeneous 
coordinate is not equal to zero. 
Now we take a coordinate of $\Pj^{n-1}$ 
around $x_0$ defined by

\begin{eqnarray*}
    \varphi:\C^{n-1} \ni 
            x&=&(x_1,\ldots,x_{m_1-1},x_{m_1+1},\ldots,x_n)\\
            &&\mapsto  
            q=[x_1,\ldots,x_{m_1-1},1,x_{m_1+1},\ldots,x_n] 
            \in \Pj^{n-1}. 
\end{eqnarray*}

\normalsize

\noindent 
Since we chose $x_0$ so that $H_p(x_0)\ne 0$, 
it holds that the derivative $\frac{\partial\det A_p}{\partial q_n}$ 
does not vanish at $x_0$. 
So by the implicit function theorem, 
there exists a function $Q_n$ 
defined on a sufficiently small neighborhood of $x_0$ such that 
\begin{eqnarray*}
    & & \lambda_p(x_1,\ldots,x_{m_1-1},1,x_{m_1+1},\ldots,x_{n-1},
              Q_n(x)) \\ 
    & & \qquad = 
        \det A_p(x_1,\ldots,x_{m_1-1},1,x_{m_1+1},\ldots,x_{n-1},
             Q_n(x))=0. 
\end{eqnarray*}
Since 
$$
    x_{m_1} = 1, \quad 
    x_{m_2} = \frac{c_{m_2}}{c_{m_1}}, \quad 
    x_{m_3} = \frac{c_{m_3}}{c_{m_1}} \qquad 
    \mbox{on} \quad \V^{n-3}, 
$$
$(x_{j_1},\ldots,x_{j_{n-4}})$ 
is considered as a local coordinate system of the variety $X^{n-4}$ 
around the regular point $x_0$. 
Since 
$$
     \frac{\partial Q_n}{\partial x_{j_\ell}} 
      = -\frac{\partial\det A_p}{\partial q_{j_\ell}} 
        /\frac{\partial\det A_p}{\partial q_n} 
        \qquad (\ell=1,\ldots,n-4) 
$$
holds, 
one can easily check that 
the condition $\det S_p(x_0)\ne 0$ implies that 
the matrix 
$$
    \left(\frac{\partial(\stackrel\circ{f_p^k}\circ\varphi)}
               {\partial x_{j_\ell}} 
         +\frac{\partial Q_n}{\partial x_{j_\ell}} 
          \frac{\partial(\stackrel\circ{f_p^k}\circ\varphi)}
               {\partial x_n} 
    \right)_{k=1,\ldots,n-1;\,\ell=1,\ldots,n-4} 
$$
is of rank $n-4$ at $x_0$. 
Hence the Jacobi matrix of $\Fl_p$ is of rank $n-4$ at $x_0$, 
and so is that of $\widehat{\Fl_p}$ at $\pi^{-1}(x_0)$. 
Thus by the proper mapping theorem, 
$\widehat{\Fl_p}(\hat X^{n-4})$ 
is an analytic subset of dimension $n-4$ 
in the same dimensional complex projective space $\W_p^{n-4}$. 
Hence $\widehat{\Fl_p}(\hat X^{n-4})=\W_p^{n-4}$. 
\end{Proof}

\medskip

\begin{Lemma} \label{lem:3} 
Let 
$\W_{\R}=\{[a]\in \W_p^{n-4} \, ; \, a_j \in \R \}$. 
Then 
$$
    \left\{\W^{n-4}_p \setminus 
           \widehat{\Fl_p}(\Z(H_p\circ \pi)\cap \hat X^{n-4})
    \right\} \cap \W_{\R} 
$$
is an open dense subset in $\W_{\R}$. 
\end{Lemma}

\begin{Proof}
By the proper mapping theorem and the theorem of Chow, 
$\widehat{\Fl_p}(\Z(H_p \circ \pi))$ 
is an algebraic subset of $\W^{n-4}_p$. 
Thus it is a closed subset in $\W^{n-4}_p$ 
in the usual topology. 
Hence 
$\left\{\W^{n-4}_p \setminus 
 \widehat{\Fl_p}(\Z(H_p\circ \pi) \cap \hat X^{n-4}) 
 \right\} \cap \W_{\R}$ 
is an open subset in $\W_{\R}$. 
Suppose that it is not dense in $\W_{\R}$. 
We may assume that  
$\widehat{\Fl_p}(\Z(H_p\circ \pi)\cap \hat X^{n-4})$ 
is common zeros of some homogeneous polynomials 
$\bigcap_{j=1}^r \Z(h_j)$. 
Then there exists an open subset in $\W_p^{n-4}$ on which 
each $h_j$ is identically zero. 
Since $\W_{\R}$ is a totally real subset 
of the complex projective space $\W_p^{n-4}$, 
by Lemma A.2 in Appendix 
we have
$$
    h_1=\cdots=h_r=0 \qquad \mbox{on} \quad \W_p^{n-4}. 
$$
This implies that 
$\widehat{\Fl_p}(\Z(H_p\circ \pi)\cap \hat X^{n-4}) = \W_p^{n-4}$. 
So it holds that 
\begin{eqnarray*}
    n-4 &  =  & \dim\W^{n-4}_p 
           =    \dim\widehat{\Fl_p}(\Z(H_p\circ \pi)\cap \hat X^{n-4}) \\ 
        & \le & \dim \Z(H_p\circ \pi)\cap \hat X^{n-4} 
          \le   \dim\hat X^{n-4} 
           =    n-4. 
\end{eqnarray*}
By the irreducibility of $\hat X^{n-4}$, 
we have 
$\Z(H_p\circ \pi)\cap \hat X^{n-4}=\hat X^{n-4}$. 
But this contradicts the fact that $H_p(x_0)\ne 0$. 
\end{Proof}

\medskip

\noindent
{(\bf Proof of Proposition \ref{Prop:A})} 
Let $p$ be a point in $U$ satisfying $\dim \W_p^{n-4}=n-4$.
As we mentioned before, 
$\dim\W_p^{n-4}=n-4$ 
holds on an open dense subset of $\{p\in\C^n\}$. 
Then for any
$$
    [a] \in \left(\W^{n-4}_p\setminus 
                  \widehat{\Fl_p}(\Z(H_p\circ \pi) \cap \hat X^{n-4}) 
            \right) \cap \W_{\R}, 
$$
there exists 
$x \in X^{n-4}\setminus \Z(H_p)$ 
such that $\Fl_p(x)=[a]$ 
by Lemma \ref{lem:2} and Lemma \ref{lem:3}. 
Since $\fl_p^j(x)\ne 0$ and 
also the resultant $R_p(x)$ does not vanish, 
$(x,b_p(x))$ induces an $n$-end catenoid 
with the flux data $(p,a)$ 
by Theorem \ref{Thm:1}. 
\hfill (q.e.d.)

\medskip

\vskip .3cm

For the later application, 
the following modification of Proposition \ref{Prop:A} 
is needed: 
Recall here that 
any elements of the matrices $A_p$ and $J_p$ are rational functions in 
$p_1,\ldots,p_n$, 
$\bar p_1,\ldots,\bar p_n$ and $q_1,\ldots,q_n$. 
Let $\check A_p$ and $\check J_p$ be the matrices 
obtained by replacing $\bar p_n$ by $p_n$, 
namely 
\begin{eqnarray}
    \label{chk1} 
    \check A_p 
     & := & A_p(p_1,\ldots,p_n,\bar p_1,\ldots,\bar p_{n-1},p_n, 
                      q_1,\ldots,q_n), \\ 
    \label{chk2} 
    \check J_p 
     & := & J_p(p_1,\ldots,p_n,\bar p_1,\ldots,\bar p_{n-1},p_n,
                q_1,\ldots,q_n), 
\end{eqnarray}
and let $\check b_p^j$ (resp. $\check f_p^j$, $\check \W_p^{n-4}$) 
be the vector (resp. function, set) 
obtained by replacing $\bar p_n$ 
in $b_p^j$ (resp. $f_p^j$, $\W_p^{n-4}$) 
by $p_n$. 
\medskip

\begin{Prop}\label{OOO}
Suppose that there exist 
$u_0\in \C^{n-1}\times \R$ and 
a point $c=[c_1,\ldots,c_n]\in \Pj^{n-1}$ 
satisfying the following conditions: 
\begin{enumerate}
\item[\rm (1)]
    $c_1,\ldots,c_n$ are all distinct{\rm ;} 
\item[{\rm (2)}] 
    The rank of the matrix $\check A_{u_0}(c)$ is $n-1${\rm ;} 
\item[\rm (3)]
    $\frac{\partial\det\check A_{u_0}}{\partial q_n}$ 
    does not vanish at $q=c${\rm ;}
\item[{\rm (4)}] 
    The rank of the matrix $\check J_{u_0}(c)$ is $n-4${\rm ;} 
\item[{\rm (5)}] 
    Two polynomials $P(z)$ and $Q(z)$ 
    defined in {\rm (\ref{Def:P})} and {\rm (\ref{Def:Q})} 
    associated with the data $(q,p)=(c,u_0)$ 
    and $\check b=\check b_{u_0}(c)$ 
    are mutually prime and 
    one of them has degree $n-1${\rm ;}
\item[{\rm (6)}] 
    $\check f_{u_0}^j(c)\ne 0$ $(j=1,\ldots,n)${\rm ;} 
\item[{\rm (7)}] 
    $c_j\ne 0$ $(j=1,\ldots,n-1)$. 
\end{enumerate}
Then there exists an open dense subset $U\subset \C^{n-1}\times \R$ 
and an open dense subset $\Omega_p$ of the totally
real set $\W_{\R}=\{[a]\in\check\W_p^{n-4}\,;\,a_j\in\R\}$
such that, 
for $p=(p_1,\ldots,p_n)\in U$
and $[a]\in\Omega_p$,
there exists an {\rm(}non-branched{\rm )} $n$-end catenoid with the flux
data $(p,a)$.
\end{Prop}

\begin{Proof}
The proof of Proposition \ref{Prop:A} works 
even if we replace $\bar p_n$ by $p_n$. 
When $p_n$ is real, 
$\check A_p$, $\check J_p$, $\check \Fl_p$ and $\check \W_p^{n-4}$ 
coincide with
$A_p$, $J_p$, $\Fl_p$ and $\W_p^{n-4}$ respectively. 
In fact, by the same proof as Lemma 1.7, we can prove that
$U:=\{p\in \C^{n-1}\times \R\,;\,
\, \Z(\lambda_p)\cap \V^{n-3}\not\subset 
    \Z(H_p)\}
$
is open dense in $\C^{n-1}\times \R$, because we only need 
the real analyticity with respect to the parameter $p$
for applying Lemma A.1.
\end{Proof}

\medskip

\begin{Rmk} \label{PPP} 
To solve the inverse problem of the flux formula, 
we may assume that $p_n \in \R$ 
since by a suitable rotation in $\{(x,y,z)\in\R^3\}$, 
we can choose that $v_n$ is in the $xz$-plane. 
\end{Rmk}

%% file: gen2.tex
\section{Finding a regular point of the flux map}

In the previous section, 
we reduced our inverse problem 
to finding a regular point of the flux map. 
However, 
the following difficulties arise in this process. 
\begin{itemize}
\item As seen in \cite{kat} and \cite{kuy}, 
      $n$-end catenoids with many symmetries 
      are easy to construct. 
      But unfortunately, 
      they are not expected to be 
      a regular point of the flux map 
      because of their symmetries. 
\item If we take a less symmetric $n$-end catenoid, 
      the computation of the rank of the flux map 
      is very complicated 
      and hard to calculate even by computer. 
\end{itemize}
To avoid these difficulties, 
we first take an $n$-end catenoid 
with many symmetries, 
and next consider a perturbation of it 
which attains the desired properties. 
 \par
Set $m:=n-1$. 
First we consider a $1$-parameter family of 
$(m+1)$-end catenoids given in \cite{kat}; 
\begin{equation} \label{sweeper} 
    \left\{\begin{array}{l}
     p_j := r\zeta^{j-1} \qquad (j=1,\ldots,m), \quad 
     p_{m+1} := 0, \\ 
     a^1 = \cdots = a^m := \displaystyle{\frac{m-1}{2}}r(r^2+1), \quad 
     a^{m+1} := \displaystyle{\frac{m(m-1)}{2}}r(r^2-1), \\ 
     q_j := \zeta^{j-1} \qquad (j=1,\ldots,m), \quad 
     q_{m+1}:=0, \\ 
     b^1 = \cdots = b^m :=1, \quad 
     b^{m+1} := \displaystyle{\frac{m-1}{2}}(r^2-1),
    \end{array}\right. 
\end{equation}
where $r>0$, 
$r\ne 1$ and $\zeta:=\exp(2\pi\sqrt{-1}/m)$. 
In fact, 
they are $(m+1)$-end catenoids without branch points 
by Remark \ref{rmk:1}, 
and are invariant under the action of the cyclic group $Z_m$. 
Unfortunately, 
as we shall see below, 
$J_p(q)=\check J_p(q)=0$ holds for any of these examples, 
namely they all are singular points of the flux maps. 
However, 
we will show that there exists a regular point near them. 
 \par
Note here that 
the matrix $A_p(q)$ (defined in (\ref{A})) 
for the example above is given by 

 \medskip

\footnotesize

\begin{equation} \label{eq:mat1} 
    A_p(q) 
     = \pmatrix{
           0 & \frac{1+r^2\zeta^1}{q_1-q_2} & \dots 
             & \frac{1+r^2\zeta^{m-1}}{q_1-q_m} 
             & \frac{1}{q_1-q_{m+1}} \cr 
           \frac{1+r^2\zeta^{-1}}{q_2-q_1} & 0 & \dots 
             & \frac{1+r^2\zeta^{m-2}}{q_2-q_m} 
             & \frac{1}{q_2-q_{m+1}} \cr 
           \vdots & \vdots & \ddots & \vdots & \vdots \cr 
           \frac{1+r^2\zeta^{-(m-1)}}{q_m-q_1} 
             & \frac{1+r^2\zeta^{-(m-2)}}{q_m-q_2} 
             & \dots & 0 
             & \frac{1}{q_m-q_{m+1}} \cr 
           \frac{1}{q_{m+1}-q_1} 
             & \frac{1}{q_{m+1}-q_2} & \dots 
             & \frac{1}{q_{m+1}-q_m} & 0 \cr 
               }. 
\end{equation}

\normalsize

 \medskip

Now, 
We consider a $1$-parameter family of matrices 

 \medskip

\footnotesize

\begin{equation} \label{eq:matS} 
    A(q,\mu) 
     := \pmatrix{
         0 & \frac{1+\mu\zeta^1}{q_1-q_2} & \dots 
           & \frac{1+\mu\zeta^{m-1}}{q_1-q_m} & \frac{1}{q_1-q_{m+1}} \cr 
         \frac{1+\mu\zeta^{-1}}{q_2-q_1} & 0 & \dots 
           & \frac{1+\mu\zeta^{m-2}}{q_2-q_m} & \frac{1}{q_2-q_{m+1}} \cr 
         \vdots & \vdots & \ddots & \vdots & \vdots \cr 
         \frac{1+\mu\zeta^{-(m-1)}}{q_m-q_1} 
           & \frac{1+\mu\zeta^{-(m-2)}}{q_m-q_2} & \dots & 0 
           & \frac{1}{q_m-q_{m+1}} \cr 
         \frac{1}{q_{m+1}-q_1} & \frac{1}{q_{m+1}-q_2} & \dots 
           & \frac{1}{q_{m+1}-q_m} & 0 \cr 
        }.
\end{equation}

\normalsize

 \medskip

By comparing (\ref{eq:mat1}) with (\ref{eq:matS}), 
we have $A(q,r^2)=A_p(q)$ for $p$ as in (\ref{sweeper}). 
When we evaluate it at 
$q=q^0:=(1,\zeta^1,\ldots,\zeta^{m-1},0)$, 
we have 

 \medskip

\footnotesize

\begin{equation} \label{eq:mat2} 
    A(q^0,\mu) 
      = \pmatrix{
         0 & \frac{1+\mu\zeta^1}{1-\zeta^1} & \dots 
           & \frac{1+\mu\zeta^{m-1}}{1-\zeta^{m-1}} & 1 \cr 
         \frac{1+\mu\zeta^{-1}}{\zeta^1-1} & 0 & \dots 
           & \frac{1+\mu\zeta^{m-2}}{\zeta^1-\zeta^{m-1}} & \zeta^{-1} \cr
         \vdots & \vdots & \ddots & \vdots & \vdots \cr 
         \frac{1+\mu\zeta^{-(m-1)}}{\zeta^{m-1}-1} 
           & \frac{1+\mu\zeta^{-(m-2)}}{\zeta^{m-1}-\zeta^1} 
           & \dots & 0 & \zeta^{-(m-1)} \cr 
         -1 & -\zeta^{-1} & \dots & -\zeta^{-(m-1)} & 0 \cr 
        }. 
\end{equation}

\normalsize

 \medskip

\noindent

We remark that 
the matrix $A(q^0,\mu)$ has the simplest form when $\mu=-1$. 
The following lemma holds. 

 \medskip

\begin{Lemma} \label{lem:ker} 
The $(m+1)$-matrix $A(q^0,\mu)$ is of rank $m$ 
except for finite values of $\mu\in \R$. 
Moreover $A(q^0,\mu)$ has a $0$-eigenvector given by 
$$
    {}^t\left(1,\ldots,1,\frac{m-1}{2}(\mu-1)\right). 
$$
\end{Lemma}

\begin{Proof} 
The second assertion is easily checked. 
Hence the rank of the matrix $A(q^0,-1)$ is at most $m$. 
Moreover, 
it is easy to see that 
the rank of the matrix $A(q^0,-1)$ is $m$. 
Since each component of $A(q^0,\mu)$ is a polynomial in $\mu$, 
the first assertion is obtained. 
\end{Proof} 

 \medskip

\begin{Rmk} \label{rmk:ker} 
Similarly, 
a $0$-eigenvector of ${}^tA(q^0,\mu)$ is given by 
$$
    {}^t\left(1,\ldots,1,\frac{1}{2}\{2\mu-(m-1)(\mu+1)\}\right). 
$$
\end{Rmk}

 \medskip

\begin{Prop} \label{prop:ess} 
The following identity holds. 
$$
    \frac{\partial\det A}{\partial q_j}(q^0,\mu)
     = 0 
    \qquad (j=1,\ldots,m+1). 
$$
\end{Prop}

\begin{Proof} 
We denote the cofactor matrix of $A(q,\mu)$ by $B(q,\mu)$.
By Lemma \ref{lem:ker} and Remark \ref{rmk:ker}, 
it can be easily checked that 
$B(q^0,\mu)$ is written in the form 
$B(q^0,\mu)=f(\mu)S(\mu)$, 
where $f(\mu)$ is a polynomial in $\mu$ 
satisfying $f(-1)=1$, 
\begin{equation} \label{B} 
    S(\mu) := \pmatrix{1 & \cdots & 1 & \psi(\mu) \cr 
                       \vdots & \ddots & \vdots & \vdots \cr 
                       1 & \cdots & 1 & \psi(\mu) \cr 
                       \varphi(\mu) & \cdots & \varphi(\mu) 
                         & \varphi(\mu)\cdot \psi(\mu) 
                      }, 
\end{equation}
and $\varphi(\mu)$ and $\psi(\mu)$ are explicitly given by
$$
    \varphi(\mu) := \frac{m-1}{2}(\mu-1), \qquad 
    \psi(\mu) := \frac{1}{2}\biggl\{2\mu-(m-1)(\mu+1)\biggr\}. 
$$
 \par
Note here that 
$$
    \frac{\partial\det A}{\partial q_j}(q,\mu) 
     = \tr\left(\frac{\partial A}{\partial q_j}(q,\mu) 
                \cdot B(q,\mu) 
          \right) 
$$
always holds for any $j$. 
Denote the $(k,\ell)$-component of the matrix $A(q,\mu)$ 
by $\alpha_{k\ell}(q,\mu)$. 
Then we have 
$$
    \frac{\partial\alpha_{k\ell}}{\partial q_j}(q^0,\mu) 
     = \left\{
       \begin{array}{ll} 
           -\displaystyle{ 
                \frac{1+\mu\zeta^{\ell-j}} 
                     {(\zeta^{j-1}-\zeta^{\ell-1})^2} 
                         } 
             & (k=j;\ell=1,\ldots,m;\ell\ne j) \\ 
           -\zeta^{-2(j-1)}  
             & (k=j;\ell=m+1) \\ 
            \displaystyle{ 
                \frac{1+\mu\zeta^{j-k}}
                     {(\zeta^{k-1}-\zeta^{j-1})^2} 
                         } 
             & (k=1,\ldots,m;k\ne j;\ell=j) \\ 
            \zeta^{-2(j-1)} 
             & (k=m+1;\ell=j) \\ 
           0 & \mbox{elsewhere} 
       \end{array}
       \right. 
$$
for $j=1,\ldots,m$, 
and 
$$ 
    \frac{\partial\alpha_{kl}}{\partial q_{m+1}}(q^0,\mu) 
     = \left\{\begin{array}{ll} 
           \zeta^{-2(k-1)}    & (k=1,\ldots,m;\ell=m+1) \\ 
          -\zeta^{-2(\ell-1)} & (k=m+1;\ell=1,\ldots,m) \\ 
           0                  & \mbox{elsewhere} 
              \end{array}
       \right. 
$$
for $j=m+1$. 
 \par 
For $j=1,\ldots,m$, 
by using the formula above, 
we have 
\begin{eqnarray*}
     &   & 
    \frac{\partial\det A}{\partial q_j}(q^0,\mu)
       =   \tr\left(\frac{\partial A}{\partial q_j}(q^0,\mu) 
                    \cdot B(q^0,\mu) 
              \right) \\ 
     &   & \quad
       =   \Condsum{k=1}{k\ne j}^{m} 
            f(\mu)\frac{\partial\alpha_{kj}}{\partial q_j}(q^0,\mu) 
          +\Condsum{\ell=1}{\ell\ne j}^{m} 
            f(\mu)\frac{\partial\alpha_{j\ell}}{\partial q_j}(q^0,\mu) \\ 
     &   & \qquad 
          +\frac{\partial\alpha_{j m+1}}{\partial q_j}(q^0,\mu) 
            f(\mu)\varphi(\mu) 
          +\frac{\partial\alpha_{m+1 j}}{\partial q_j}(q^0,\mu) 
            f(\mu)\psi(\mu) \\ 
     &   & \quad 
       =   f(\mu)\left\{
           \Condsum{k=1}{k\ne j}^{m} 
            \frac{1+\mu\zeta^{j-k}}{(\zeta^{k-1}-\zeta^{j-1})^2} 
          -\Condsum{\ell=1}{\ell\ne j}^{m} 
            \frac{1+\mu\zeta^{\ell-j}}{(\zeta^{j-1}-\zeta^{\ell-1})^2} 
                                   -\zeta^{-2(j-1)}\varphi(\mu) 
           +\zeta^{-2(j-1)}\psi(\mu) 
           \right\} \\ 
     &   & \quad 
       =   f(\mu)\zeta^{-2(j-1)}
            \left\{
             \sum_{k=1}^{m-1} 
              \frac{1+\mu\zeta^{-k}}{(1-\zeta^k)^2} 
            -\sum_{k=1}^{m-1} 
              \frac{1+\mu\zeta^k}{(1-\zeta^k)^2} 
            -(m-2)\mu 
            \right\} \\ 
     &   & \quad 
       =   \mu f(\mu)\zeta^{-2(j-1)}
            \left\{
             \sum_{k=1}^{m-1} 
                \frac{1+\zeta^k}{\zeta^k(1-\zeta^k)} 
             -(m-2) 
            \right\} \\ 
     &   & \quad 
       =   \mu f(\mu)\zeta^{-2(j-1)}
            \left\{
              \sum_{k=1}^{m-1} 
                \frac{1}{\zeta^k} 
             +\sum_{k=1}^{m-1} 
                \frac{2}{1-\zeta^k} 
             -(m-2) 
            \right\} \\ 
     &   & \quad 
       =   \mu f(\mu)\zeta^{-2(j-1)}
            \left\{-1+(m-1)-(m-2) 
            \right\} = 0. 
\end{eqnarray*}
On the other hand, 
for $j=m+1$, 
we have 
\begin{eqnarray*}
    \frac{\partial\det A}{\partial q_{m+1}}(q^0,\mu) 
     & = & \tr\left(\frac{\partial A}{\partial q_{m+1}}(q^0,\mu) 
                    \cdot B(q^0,\mu) 
              \right) \\ 
     & = & \sum_{k=1}^{m} \zeta^{-2(k-1)}f(\mu)\varphi(\mu) 
          -\sum_{\ell=1}^{m} \zeta^{-2(\ell-1)}f(\mu)\psi(\mu) \\ 
     & = & f(\mu)(\varphi(\mu)-\psi(\mu)) 
          \sum_{k=1}^{m} \zeta^{-2k} =0. 
\end{eqnarray*}
This completes the proof. 
\end{Proof}

 \medskip

By Lemma \ref{lem:ker} and Proposition \ref{prop:ess}, 
it follows that $J_{rq^0}(q^0)=0$ \,\, $(r\in \R)$. 
Therefore, 
we try to perturb a sampling point. 
To do this, 
we consider an $m$-matrix $\Gamma_{m+1}(\mu)$ by 

\footnotesize

$$
    \Gamma_{m+1}(\mu) 
      := 
        \left( 
        \frac{\partial^2 \det A}{\partial q_1\partial q_{m+1}}(q^0,\mu) 
        \cdot 
        \frac{\partial(f^k/f^{m+1})}{\partial q_j}(q^0,\mu) 
       -\frac{\partial^2 \det A}{\partial q_1\partial q_j}(q^0,\mu) 
        \cdot 
        \frac{\partial(f^k/f^{m+1})}{\partial q_{m+1}}(q^0,\mu) 
        \right)_{j,k=1,\ldots,m}, 
$$

\normalsize

\noindent 
where we denote 
the $(j,k)$-component of the cofactor matrix $B(q,\mu)$ 
by $\beta_{jk}(q,\mu)$, 
and set 

\footnotesize

\begin{eqnarray} \label{fj} 
     &    & \\ 
    f^k(q,\mu) 
     & := & \beta_{k\,m+1}(q,\mu) 
            \left( 
            \Condsum{j=1}{j\ne k}^m \beta_{j\,m+1}(q,\mu) 
            \frac{\zeta^{k-1}-\zeta^{j-1}}{q_k-q_j} 
           +\beta_{m+1\,m+1}(q,\mu) 
            \frac{\zeta^{k-1}}{q_k-q_{m+1}} 
            \right) 
            \nonumber \\ 
     &    & \qquad (k=1,\ldots,m), 
            \nonumber \\ 
    f^{m+1}(q,\mu) 
     & := & \beta_{m+1\,m+1}(q,\mu) 
            \sum_{j=1}^m \beta_{j\,m+1}(q,\mu) 
            \frac{-\zeta^{j-1}}{q_{m+1}-q_j}. 
            \nonumber 
\end{eqnarray}

\normalsize

\noindent 
(Compare with the definition of the matrix $J_p(q)$ and $f_p^k(q)$.)
We prove the following 

 \medskip

\begin{Thm} \label{Main} 
Suppose that 
there exists a positive number $\mu$ such that 
the matrix $\Gamma_{m+1}(\mu)$ $(n=m+1\ge 5)$ 
is of rank $m-3(=n-4)$ and
\begin{equation} \label{27}
    \frac{\partial^2 \det A}{\partial q_1\partial q_{m+1}}(q^0,\mu) 
     \ne 0. 
\end{equation}
Then, 
for each of almost all flux data, 
there exists an $n$-end catenoid 
with the flux data. 
\end{Thm}

 \medskip

Till now, 
we fix the parameter $p_{m+1}$ at 
$$
    p_{m+1}=0.
$$
Let us now move $p_{m+1}$ as a complex parameter. 

 \medskip

\begin{Lemma} \label{NNN} 
Let $\mu\ne 1$ be a positive real number 
such that $f(\mu)\ne 0$, 
where $f(\mu)$ is a polynomial given by (\ref{B}). 
Then 
$$
    \frac{\partial  \det \check A_p(q)}{\partial p_{m+1}} 
     \ne 0 
$$
at the point 
$q=q^0=(1,\zeta^1,\ldots,\zeta^{m-1},0)$ 
for $p=\sqrt{\mu}q^0$, where $\check A_p(q)$ is defined in $(\ref{chk1})$.
\end{Lemma} 

\begin{Proof} 
We denote the cofactor matrix of $\check A_p(q)$ by $\check B_p(q)$. 
Since $\check A_{\sqrt{\mu}q^0}(q)=A_{\sqrt{\mu}q^0}(q)$ 
for any $\mu>0$, 
it holds that 
$\check B_{\sqrt{\mu}q^0}(q)=B_{\sqrt{\mu}q^0}(q)$ 
and in particular, 
we have 
$\check B_{\sqrt{\mu}q^0}(q^0)=B_{\sqrt{\mu}q^0}(q^0)=B(q^0,\mu)$. 
Then we have 
$$
    \left. 
    \frac{\partial \det\check A_p(q^0)}{\partial p_{m+1}} 
    \right|_{p={\sqrt{\mu}q^0}} 
     = \tr\left( 
                  \left. 
                  \frac{\partial \check A_p(q^0)}{\partial p_{m+1}} 
                  \right|_{p={\sqrt{\mu}q^0}} 
                   \cdot 
                  B_{\sqrt{\mu}q^0}(q^0) 
          \right). 
$$
Since 

\footnotesize

$$
    \mbox{the $(j,k)$-component of } 
    \left. 
    \frac{\partial \check A_p(q^0)}{\partial p_{m+1}}
    \right|_{p={\sqrt{\mu}q^0}} 
     = \left\{ 
       \begin{array}{ll} 
            \zeta^{-2(j-1)} & (j=1,\ldots,m;k=m+1) \\ 
            -1              & (j=m+1;k=1,\ldots,m) \\ 
            0               & \mbox{elsewhere}, 
       \end{array} 
       \right. 
$$

\normalsize

\noindent 
by (\ref{B}), 
we have 
\begin{eqnarray*}
    & & \tr\left( 
           \left. 
             \frac{\partial  \check A_p(q^0)}{\partial p_{m+1}} 
           \right|_{p={\sqrt{\mu}q^0}} 
             \cdot B(q^0,\mu) 
           \right) \\ 
    & & \qquad 
        = f(\mu)\left\{ 
             \varphi(\mu)\sum_{k=1}^{m}\zeta^{-2(k-1)} 
             -(m-1)\psi(\mu) 
             \right\} \\ 
    & & \qquad 
        = -(m-1)f(\mu)\psi(\mu) 
        = \frac{(m-1)^2}{2}(\mu-1)f(\mu) \ne 0. 
\end{eqnarray*}
Now the assertion is clear. 
\end{Proof}

 \medskip

\noindent 
{\bf (Proof of Theorem \ref{Main}.)} 
Since $f(\mu)$ is a polynomial in $\mu$ and $f(\mu)\not\equiv 0$, 
by our assumptions and Lemmas \ref{lem:ker} and \ref{NNN}, 
we can choose a positive number $\mu$ such that 
$f(\mu)\ne 0$, 
rank $\check A_{\sqrt{\mu}q^0}(q^0)=m$, 
rank $\Gamma_{m+1}(\mu)=m-3$, 
$$
    \frac{\partial^2\det \check A_{\sqrt{\mu}q^0}}
         {\partial q_1\partial q_{m+1}}(q^0) \ne 0 
     \qquad\mbox{and}\qquad 
    \left. 
    \frac{\partial \det\check A_p(q^0)}{\partial p_{m+1}} 
    \right|_{p={\sqrt{\mu}q^0}} 
     \ne 0. 
$$
Throughout this proof, 
we fix the parameters except for $q_1$ and $p_{m+1}$ 
to the same values as $q=q^0$ and $p=\sqrt{\mu}q^0$: 
\begin{eqnarray*}
    & & p_j = \sqrt{\mu}\zeta^{j-1}\qquad (j=1,\ldots,m), \\ 
    & & q_j = \zeta^{j-1} \qquad (j=2,\ldots,m), 
        \qquad q_{m+1} = 0. 
\end{eqnarray*}
Regard $\det \check A_p(q)$ 
as a function with respect to 
only $q_1$ and $p_{m+1}$, 
and apply the implicit function theorem to the point 
$(q_1,p_{m+1})=(1,0)$. 
Then there exist 
an open neighborhood $U\subset \C$ of $1\in \C$ 
and 
a complex analytic function $p_{m+1}=p_{m+1}(q_1):U\to \C$ 
such that 
$p_{m+1}(1)=0$ 
and 
$$
    \biggl.\det \check A_p\biggr|_{p_{m+1}=p_{m+1}(q_1)} = 0 
     \qquad (q_1\in U). 
$$
Since rank $\check A_{\sqrt{\mu}q^0}(q^0)=m$, 
rank $\check A_p|_{p_{m+1}=p_{m+1}(q_1)}=m$ 
holds also for $q_1$ near $1$. 
 \par
Since $\hat A=A$ at $p=\sqrt{\mu}q^0$, by Lemma \ref{prop:ess},
we have 
$$
    \frac{\partial\det \check A_{\sqrt{\mu}q^0}}
         {\partial q_j}(q^0) 
     = 0 
     \qquad (j=1,\ldots,m+1). 
$$
On the other hand, 
the assumption (\ref{27}) yields 
$$
    \biggl. 
    \frac{\partial\det\check A_p}{\partial q_{m+1}} 
    \biggr|_{p_{m+1}=p_{m+1}(q_1)} 
     \ne 0 
$$
for any $q_1\ne 1$ enough close to $1$. 
Therefore we have 
\begin{eqnarray*}
    & & \lim_{q_1\to 1} 
        \left( 
        \left. 
        \left\{ 
            \frac{\partial(\check f_p^k/\check f_p^{m+1})} 
                 {\partial q_j} 
             - 
            \frac{\frac{\partial\det \check A_p}{\partial q_j}} 
                 {\frac{\partial\det \check A_p}{\partial q_{m+1}}} 
             \cdot 
            \frac{\partial(\check f_p^k/\check f_p^{m+1})} 
                 {\partial q_{m+1}} 
        \right\} 
        \right|_{p_{m+1}=p_{m+1}(q_1)} 
        \right)_{j,k=1,\ldots,m} 
\\ 
    & & \qquad = 
        \left( 
        \left. 
        \left\{ 
            \frac{\partial(\check f_p^k/\check f_p^{m+1})} 
                 {\partial q_j} 
             - 
            \frac{\frac{\partial^2 \det \check A_p}
                       {\partial q_1\partial q_j}} 
                 {\frac{\partial^2 \det \check A_p}
                       {\partial q_1\partial q_{m+1}}} 
             \cdot 
            \frac{\partial(\check f_p^k/\check f_p^{m+1})} 
                 {\partial q_{m+1}} 
        \right\} 
        \right|_{p=\sqrt{\mu}q^0;q=q^0} 
        \right)_{j,k=1,\ldots,m} 
         \\ 
    & & \qquad = 
        \left(
        \frac{\partial^2 \det \check A_{\sqrt{\mu}q^0}} 
             {\partial q_1\partial q_{m+1}}
             (q^0) 
        \right)^{-1} 
        \Gamma_{m+1}(\mu), 
\end{eqnarray*}
and hence 
\begin{eqnarray*}
     & & \left. 
         \mbox{rank }\check J_p 
         \right|_{p_{m+1}=p_{m+1}(q_1)} 
          \\ 
     & & \quad = 
         \mbox{rank }
         \left( 
         \left. 
         \left\{ 
         \frac{\partial(\check f_p^k/\check f_p^{m+1})} 
              {\partial q_j} 
          - 
         \frac{\frac{\partial\det \check A_p}{\partial q_j}} 
              {\frac{\partial\det \check A_p}{\partial q_{m+1}}} 
          \cdot 
         \frac{\partial(\check f_p^k/\check f_p^{m+1})} 
              {\partial q_{m+1}} 
         \right\} 
         \right|_{p_{m+1}=p_{m+1}(q_1)} 
         \right)_{j,k=1,\ldots,m} 
          \\ 
    & & \quad = m-3 = n-4 
\end{eqnarray*}
for any $q_1$ as above. 
 \par
Since the initial sampling point 
$q=q^0$, 
$p=\sqrt{\mu}q^0$ 
is chosen from the data 
which realizes a non-branched $n$-end catenoid $(n=m+1)$, 
$\Delta(q^0)\ne 0$ and $q^0_j\ne 0$ $(j=1,\ldots,m)$, 
the other conditions in Proposition \ref{OOO} 
are also satisfied for $q_1$ near $1$ such that $p_{m+1}\in \R$. 
Now, 
by Remark \ref{PPP}, 
we have proved the theorem.
\hfill(q.e.d.)

 \medskip

Thus we will get our main theorem in Introduction, 
if the matrix $\Gamma_{m+1}(\mu)$ is of rank $m-3(=n-4)$ 
and (\ref{27}) holds for some $\mu>0$, 
which will be shown in the next section.

%% file: gen3.tex
\section{Computation of $\Gamma_{m+1}(\mu)$}

In this section, 
we compute the matrix $\Gamma_{m+1}(\mu)$ 
defined in the previous section, 
and show that it is of rank $m-3$ for almost all $\mu>0,\ne 1$. 

 \medskip

\paragraph*{(Computation of 
$\frac{\partial f^k}{\partial q_j}(q^0,\mu)$)} 
As before, 
we write $A(q,\mu)=:(\alpha_{k\ell})_{k,\ell=1,\ldots,m+1}$ 
and $B(q,\mu)=:(\beta_{k\ell})_{k,\ell=1,\ldots,m+1}$. 
By (\ref{fj}), 
(\ref{B}) 
and straightforward calculations, 
we have, 
for any $k=1,\ldots,m$, 
\begin{equation} \label{AS2-} 
    \frac{\partial f^k}{\partial q_j} 
       =   f\psi
           \left[ 
               (m-1+\varphi) 
               \frac{\partial\beta_{k\, m+1}}{\partial q_j} 
                + 
               \Condsum{\ell=1}{\ell\ne k}^{m+1} 
               \frac{\partial\beta_{\ell\,\, m+1}}{\partial q_j} 
                + 
               f\psi\zeta^{1-j}\eta_1 
           \right] 
\end{equation}
at $(q^0,\mu)$, 
where 
$$
    \eta_1(\mu) 
        := \left\{
           \begin{array}{ll}
               -\frac{m-1}{2}-\varphi(\mu) & (j=k) \\ 
               \frac{1}{\zeta^{k-j}-1}     & (j=1,\ldots,m;j\ne k) \\ 
               \zeta^{j-k}\varphi(\mu)     & (j=m+1), 
           \end{array} 
           \right. 
$$
and for $k=m+1$, 

\footnotesize

\begin{equation} \label{AS2--} 
    \frac{\partial f^{m+1}}{\partial q_j} 
       =   f\psi 
           \left[ 
               m 
               \frac{\partial\beta_{m+1\, m+1}}{\partial q_j} 
                + 
               \varphi 
               \sum_{\ell=1}^{m} 
               \frac{\partial\beta_{\ell\,\, m+1}}{\partial q_j} 
                - 
               \left\{ 
               \begin{array}{ll} 
                   f\varphi\psi\zeta^{1-j} & (j=1,\ldots,m) \\ 
                   0                       & (j=m+1) 
               \end{array} 
               \right. 
           \right]. 
\end{equation}

\normalsize

Hence we have only to compute $f(\mu)$ and 
$\frac{\partial\beta_{k\,m+1}}{\partial q_j}(q^0,\mu)$. 
Denote the first $m\times m$-submatrix of $A(q^0,\mu)$ by $A^0(\mu)$. 
Clearly $f\varphi\psi=\beta_{m+1 m+1}=\det A^0$. 
Set $C_1:=$diag$[1,\zeta^1,\ldots,\zeta^{m-1}]$. 
Since $C_1A^0$ is a cyclic matrix 
whose $(j,k)$-component is equal to 
$(1+\mu\zeta^{k-j})/(1-\zeta^{k-j})$, 
and whose diagonal components vanish, 
it can be diagonalized as 
${C_2}^{-1}C_1A^0C_2=$diag$[\psi_1,\ldots,\psi_m]$, 
where 
$$
    C_2 := \left( 
           \begin{array}{llll} 
            1           & 1              & \cdots & 1 \\ 
            \zeta^1     & \zeta^2        & \cdots & 1 \\ 
            \vdots      & \vdots         & \ddots & \vdots \\ 
            \zeta^{m-1} & \zeta^{2(m-1)} & \cdots & 1 
           \end{array} 
           \right) 
$$
and the eigenvalues 
$\psi_1,\ldots,\psi_m$ of $C_1A^0$ 
are given by 
\begin{eqnarray*}
    \psi_\ell(\mu) 
     & = & \sum_{k=2}^m 
           \frac{1+\mu\zeta^{k-1}}{1-\zeta^{k-1}}(\zeta^l)^{k-1} \\ 
     & = & \left\{ 
           \begin{array}{ll} 
            \left(\ell-\frac{m-1}{2}\right)\mu 
             + \left(\ell-\frac{m+1}{2}\right) 
             & (\ell=1,\ldots,m-1) \\ 
            -\frac{m-1}{2}\mu+\frac{m-1}{2}
             & (\ell=m). 
           \end{array} 
           \right. 
\end{eqnarray*}
Now we have 
$$
    f\varphi\psi
     = (-1)^{m-1}\prod_{\ell=1}^m\psi_\ell. 
$$
Note here that $\psi_1=\psi$ and $\psi_m=-\varphi$ 
and that $\psi_\ell(\mu)\ne 0$ holds 
for any $\mu>0,\ne 1$ $(\ell=1,\ldots,m)$. 
 \par
To compute the derivatives 
$\frac{\partial B}{\partial q_j}(q^0,\mu)$ 
of the cofactor matrix $B(q,\mu)$, 
we apply the formula (\ref{b2}) in Appendix B 
by putting $X:=E_{m+1}$, 
where $E_{m+1}$ is the $(m+1)$-matrix
given by
$$
    E_{m+1}:=\pmatrix{0      & \cdots & 0      & 0      \cr 
                      \vdots & \ddots & \vdots & \vdots \cr 
                      0      & \cdots & 0      & 0      \cr 
                      0      & \cdots & 0      & 1 
                     }.
$$ 
For $A_t(q,\mu)=A(q,\mu)+tE_{m+1}$, 
we have already shown that 
$$
    \det A(q^0,\mu) 
     = 
    \frac{\partial\det A}{\partial q_j} 
    (q^0,\mu) 
     = 0 
$$
in Lemma \ref{lem:ker} and Proposition \ref{prop:ess}. 
Moreover we have 
$$
    \tr(E_{m+1} \cdot B(q^0,\mu)) = f(\mu)\varphi(\mu)\psi(\mu) \ne 0. 
$$
Thus we may apply (\ref{b2}), 
and get the following identity 
\begin{eqnarray} \label{bqj} 
     \frac{\partial B}{\partial q_j} 
      & = & \frac{1}{f\varphi\psi} 
            \left\{ 
              \tr\left(
               \frac{\partial A}{\partial q_j} 
                \cdot 
               \left.\frac{\partial Y_t}{\partial t}\right|_{t=0} 
                 \right) 
               \cdot B 
            \right. \\ 
      &   & \qquad 
            \left. 
             -\left.\frac{\partial Y_t}{\partial t}\right|_{t=0} 
               \cdot \frac{\partial A}{\partial q_j} 
               \cdot B 
             -B \cdot 
               \frac{\partial A}{\partial q_j} \cdot 
               \left.\frac{\partial Y_t}{\partial t}\right|_{t=0} 
            \right\} 
            \nonumber
\end{eqnarray}
at $(q^0,\mu)$, 
where $Y_t(\mu)$ is the cofactor matrix of $A(q^0,\mu)+tE_{m+1}$. 
The first $m\times m$-components of 
$\frac{\partial}{\partial t}|_{t=0}Y_t(\mu)$ 
is given as the cofactor matrix of 
the first $m\times m$-components of $A(q^0,\mu)$, 
that is 
\begin{eqnarray*}
    \det A^0\cdot (A^0)^{-1} 
     & =  & f\varphi\psi\cdot 
            C_2\mbox{diag}[{\psi_1}^{-1},\ldots,{\psi_m}^{-1}]
            {C_2}^{-1}C_1 \\ 
     & =  & \frac{f\varphi\psi}{m} 
            \left( 
            \zeta^{k-1}\sum_{\ell=1}^m\zeta^{(j-k)\ell}\psi_\ell^{-1} 
            \right)_{j,k=1,\ldots,m} \\ 
     & =: & \frac{f\varphi\psi}{m}Y^0, 
\end{eqnarray*}
and the other components of 
$\frac{\partial}{\partial t}|_{t=0}Y_t(\mu)$ 
vanish. 
Namely 
$$
    \left.\frac{\partial}{\partial t}\right|_{t=0}Y_t(\mu) 
     = \pmatrix{
           &        &   & 0      \cr 
           & \displaystyle{\frac{f\varphi\psi}{m}Y^0} 
                    &   & \vdots \cr 
           &        &   & 0      \cr 
         0 & \cdots & 0 & 0 
               }. 
$$
Therefore we have 
\begin{equation} \label{bnqj} 
     \left( 
     \frac{\partial\beta_{k\, m+1}}{\partial q_j} 
     \right)_{k=1,\ldots,m+1} 
      =     \frac{f\psi}{m} 
            \left\{ 
              \tr\left( 
               \frac{\partial A}{\partial q_j} 
                \cdot Y^0 
                 \right) 
               \cdot I 
              -Y^0 
               \cdot \frac{\partial A}{\partial q_j} 
            \right\} 
            \left( 
            \begin{array}{l} 
             1 \\ \vdots \\ 1 \\ \varphi 
            \end{array} 
            \right) 
            \nonumber
\end{equation}
at $(q^0,\mu)$. 
Recall here the values of 
$\frac{\partial \alpha_{k\ell}}{\partial q_j}(q^0,\mu)$ 
computed in the proof of Proposition \ref{prop:ess}. 
Now, 
by direct computation, 
we have 

\footnotesize

\begin{eqnarray} \label{bbqj} 
    \frac{\partial \beta_{k\, m+1}}{\partial q_j}(q^0,\mu) 
      & = & 
     -f(\mu)\psi(\mu)\zeta^{1-j} \\ 
      &   & \times 
     \left\{
     \begin{array}{ll} 
         \left(1-\frac{1}{2m}\eta_2(\mu)\right) 
           & (k,j=1,\ldots,m) \\ 
         \varphi(\mu) 
           & (k=m+1;j=1,\ldots,m) \\ 
         \zeta^{1-k}\varphi(\mu)\psi_{m-1}(\mu)^{-1} 
           & (k=1,\ldots,m;j=m+1) \\ 
         0 & (k=j=m+1), 
     \end{array} 
     \right. 
     \nonumber 
\end{eqnarray}

\normalsize

\noindent
where 
$$
    \eta_2(\mu) 
           := \left\{ 
              \begin{array}{ll} 
                  \frac{m(m-1)}{2} 
                 +\frac{\psi_1(\mu)}{\mu+1} 
                  \left\{m-1+(m+\varphi(\mu)) 
                         \sum_{\ell=1}^{m-1}\psi_\ell(\mu)^{-1} 
                  \right\}
                   & (k=j) \\ 
                  \frac{m}{\zeta^{k-j}-1} 
                 +\frac{\psi_1(\mu)}{\mu+1} 
                  \left\{-1+(m+\varphi(\mu)) 
                         \sum_{\ell=1}^{m-1}\zeta^{(k-j)\ell} 
                                            \psi_\ell(\mu)^{-1} 
                  \right\} 
                   & (k\ne j). 
              \end{array} 
              \right. 
$$
Putting it into (\ref{AS2-}) and (\ref{AS2--}), 
we get 

\footnotesize

\begin{eqnarray} \label{AS2} 
     & & 
    \frac{\partial f^k}{\partial q_j}(q^0,\mu) 
     = 
    -f(\mu)^2\psi(\mu)^2\zeta^{1-j} \\ 
     & & \qquad \times 
    \left\{
    \begin{array}{ll}
        2(m-1+\varphi(\mu))-\frac{m-2+\varphi(\mu)}{2m}\eta_2(\mu)
         -\eta_1(\mu)      & (k,j=1,\ldots,m) \\ 
        (2m+1)\varphi(\mu) & (k=m+1;j=1,\ldots,m) \\ 
        0                  & (k=1,\ldots,m+1;j=m+1). 
    \end{array} 
    \right. 
    \nonumber 
\end{eqnarray}

\normalsize

In particular, 
we have 
\begin{eqnarray*}
    \Gamma_{m+1}(\mu) 
     & = & (f^{m+1})^{-2}
           \frac{\partial^2\det A}{\partial q_1\partial q_{m+1}} 
            \cdot 
           \left( 
               f^{m+1}\frac{\partial f^k}{\partial q_j} - 
               f^k\frac{\partial f^{m+1}}{\partial q_j} 
           \right)_{k,j=1,\ldots,m} 
\end{eqnarray*}
at $(q^0,\mu)$. 
 \par

 \medskip

\paragraph*{(Computation of 
$\frac{\partial^2\det A}{\partial q_1\partial q_{m+1}}(q^0,\mu)$)} 
First we compute 

\footnotesize

$$
    \frac{\partial^2\det A}{\partial q_1\partial q_{m+1}} 
    (q^0,-1) 
      = 
    \tr\left( 
    \frac{\partial^2 A}{\partial q_1\partial q_{m+1}}(q^0,-1) 
    \cdot B(q^0,-1) 
     + 
    \frac{\partial A}{\partial q_1}(q^0,-1) 
    \cdot 
    \frac{\partial B}{\partial q_{m+1}}(q^0,-1) 
       \right). 
$$

\normalsize

\noindent 
It is easy to see that, 
\begin{eqnarray*}
    \frac{\partial^2 \alpha_{k\ell}}
         {\partial q_1\partial q_{m+1}}(q^0,-1) 
     & = & 
     \left\{
     \begin{array}{ll} 
         -2 & (k=1;\ell=m+1) \\ 
          2 & (k=m+1;\ell=1) \\ 
          0 & \mbox{elsewhere}. 
     \end{array} 
     \right. 
\end{eqnarray*}
On the other hand, 
we have 
$$
    \left.\frac{\partial}{\partial t}\right|_{t=0}Y_t(-1) 
     = \pmatrix{
         2-m    & \zeta^1      & \cdots & \zeta^{m-1} & 0      \cr 
         1      & (2-m)\zeta^1 & \cdots & \zeta^{m-1} & 0      \cr 
         \vdots & \vdots       & \ddots & \vdots      & \vdots \cr 
         1      & \zeta^1      & \cdots & (2-m)\zeta^{m-1} & 0 \cr 
         0      & 0            & \cdots & 0           & 0 
               }. 
$$
By putting these values into (\ref{bqj}), 
we have 
\begin{equation} \label{ASn} 
    \frac{\partial\beta_{k\ell}}
         {\partial q_{m+1}}(q^0,-1) 
      = 
     \left\{
     \begin{array}{ll} 
         -(m-1)\zeta^{1-k}+\zeta^{1-\ell} & (k,\ell=1,\ldots,m) \\ 
          (m-1)\zeta^{1-k}                & (k=1,\ldots,m;\ell=m+1) \\ 
         -(m-1)\zeta^{1-\ell}             & (k=m+1;\ell=1,\ldots,m) \\ 
         0                                & (k=\ell=m+1). 
     \end{array} 
     \right. 
\end{equation}
Now, 
by a straightforward calculation, 
we have 
\begin{equation} \label{dAqq} 
    \frac{\partial^2\det A} 
         {\partial q_1\partial q_{m+1}} 
    (q^0,-1) 
     = m(m-1) 
     \ne 0. 
\end{equation}
Since 
$\frac{\partial^2\det A}{\partial q_1\partial q_{m+1}}(q^0,\mu)$ 
is a polynomial in $\mu$, 
it does not vanish for any $\mu$ 
except for finite values. 
 \par

 \medskip

\paragraph*{(Computation of the rank of $\Gamma_{m+1}(\mu)$)} 

For any $\mu>0,\ne 1$ 
such that 
$\frac{\partial^2\det A}{\partial q_1\partial q_{m+1}}(q^0,\mu)\ne 0$, 
define a cyclic matrix 
$$
    \Gamma_{m+1}^0 
     := -\frac{1}{(f\psi)^2} 
        \left( 
            \frac{\partial f^k}{\partial q_j} - 
            \frac{f^k}{f^{m+1}}
            \frac{\partial f^{m+1}}{\partial q_j} 
        \right)_{k,j=1,\ldots,m} 
        \cdot C_1. 
$$
Then it is clear that 
the rank of $\Gamma_{m+1}$ 
is equal to the rank of $\Gamma_{m+1}^0$. 
The $(k,j)$-component $\gamma_{kj}$ of $\Gamma_{m+1}^0$ 
is given by 
$$
    \gamma_{kj} 
     = -\frac{m-1+\varphi}{m}
       -\frac{m-2+\varphi}{2m}\eta_2-\eta_1, 
$$
and the eigenvalues 
$\chi_1,\ldots,\chi_m$ of $\Gamma_{m+1}^0$ 
are given by 
\begin{eqnarray*}
    \chi_\ell(\mu) 
     & = & \sum_{j=1}^m\gamma_{1j}(\mu)(\zeta^\ell)^{j-1} \\ 
     & = & \left\{ 
           \begin{array}{ll} 
               -\frac{(\mu+1)\{(m-1)\mu+m+1\}(\ell-1)(\ell-m+1)} 
                     {4\psi_\ell(\mu)} 
                 & (\ell=1,\ldots,m-1) \\
               0 & (\ell=m). 
           \end{array} 
           \right. 
\end{eqnarray*}
Now it is clear that 
$\chi_\ell(\mu)\ne 0$ for $\ell=2,\ldots,m-2$, 
and $\Gamma_{m+1}^0$ is of rank $m-3$. 
Consequently, 
$\Gamma_{m+1}$ is of rank $m-3$ 
for any $\mu>0,\ne 1$ 
except for finite values. 
 \par

 \medskip

Now, 
by Theorem \ref{Main}, 
we get the following theorem: 

 \medskip

\begin{Thm}
For almost all given unit vectors 
$v=\{v_1,\ldots,v_n\}$ $(n\ge 5)$ in $\R^3$, 
and nonzero real numbers 
$a=\{a^1,\ldots,a^n\}$ 
satisfying $\sum_{j=1}^n a^j v_j=0$, 
there is a {\rm (}non-branched{\rm )} $n$-end catenoid 
$x:\C\setminus \{q_1,\ldots,q_n\}\to {\R}^3$ 
such that $\nu(q_j)=v_j$ 
and $a_j$ is the weight at the end $q_j$.
\end{Thm}

This theorem and the results for $n\le 4$ (\cite{l}, \cite{kuy})
imply our main theorem in Introduction.

%% file: gena.tex
\section*{Appendix A}

\renewcommand{\thesection}{\Alph{section}.}
\setcounter{section}{1}
\setcounter{equation}{0}
\setcounter{Thm}{0}

In this appendix, 
we give two lemmas 
on real analytic families of algebraic equations 
which are 
applied in the proof of Proposition \ref{Prop:A}. 

 \medskip

\begin{Lemma}
Let 
$\{f_p(q_1,\ldots,q_n)\}_{p\in \R^\ell}$ 
and 
$\{g_p(q_1,\ldots,q_n)\}_{p\in \R^\ell}$ 
be two real analytic families of polynomials on $\C$ 
of degree bounded by $m$. 
Suppose that 
there exists a non-empty open subset $U$ such that 
\begin{equation} \label{A1} 
    Z(f_p)\subset Z(g_p)\qquad (p\in U).
\end{equation}
Then {\rm (\ref{A1})} holds 
for all $p \in \R^\ell$ 
such that $f_p\not\equiv 0$. 
\end{Lemma}

\begin{Proof}
For each $p\in \R^\ell$, 
since the degree of $f_p$ is bounded by $m$, 
$Z(f_p)\subset Z(g_p)$ 
if and only if 
$(g_p)^m$ is divided by $f_p$. 
We operate a differential operator
$$
    D^\alpha := \frac{\partial^{|\alpha|}}
                     {\partial q_1^{\alpha_1}\cdots\partial q_n^{\alpha_n}} 
$$
into the rational function $\varphi_p:=(g_p)^m/f_p$. 
Let ${\cal N}^\alpha(\varphi_p)$ be 
a polynomial formally defined as 
$$
    {\cal N}^\alpha(\varphi_p) 
     := (f_p)^{|\alpha|+1}\cdot D^\alpha\varphi 
$$
which is the numerator part of $D^\alpha\varphi$. 
 \par
Now we fix an element $p_0 \in \R^\ell$ 
such that $f_p\not\equiv 0$, 
and choose an element $q_0 \in \C^n$ 
such that $f_{p_0}(q_0)\ne 0$. 
Since $f_p$ is real analytic 
with respect to the parameter $p$, 
we can take a subdomain $V$ of $U$ 
such that $f_p(q_0)\ne 0$ for all $p\in V$, 
and $\varphi_p$ is a polynomial on $\C$ 
of degree bounded by $m^2$ for any $p\in V$. 
Hence for any multi-index $|\alpha|>m^2$, 
we have 
${\cal N}^\alpha(\varphi_p)(q_0) = 0$ 
for $p \in V$. 
By the real analyticity 
with respect to the parameter $p$, 
we have 
${\cal N}^\alpha(\varphi_{p_0})(q_0)= 0$ 
for $|\alpha|>m^2$. 
Since $f_{p_0}(q_0)\ne 0$, 
we get 
$D^{\alpha}\varphi(q_0)=0$ 
for $|\alpha|>m^2$. 
Thus $\varphi_{p_0}$ is also a polynomial on $\C$. 
\end{Proof}

 \medskip

The following lemma is easily proved 
by using the Cauchy-Riemann equation. 

 \medskip

\begin{Lemma}
Let $\W_0$ be a totally real subset of $\Pj^{n-1}$ defined by
$$
    \W_0 := \{[a^1,\ldots,a^n]\in \Pj^{n-1}\,;\, 
               a^j\in \R \,\,(j=1,\ldots,n)\}. 
$$
Let $h$ be a homogeneous polynomial on $\C$. 
If $h$ is identically zero 
on a non-empty open subset in $\W_0$, 
then $h\equiv 0$ on $\Pj^{n-1}$. 
\end{Lemma}

%% file: genb.tex
\section*{Appendix B}

\setcounter{section}{2}
\setcounter{equation}{0}
\setcounter{Thm}{0}

Let $A$ be an $n\times n$ matrix. 
The cofactor matrix $B$ of $A$ 
is the matrix satisfying the identity 
$BA=AB=\det A\cdot I$. 
In this appendix, 
we give an identity 
which is useful to compute 
a differential of the cofactor matrix of a singular matrix. 
 \par
Let $\Omega$ be a domain in $\C$ containing the origin, 
and $A(q):\Omega \to M(n,\C)$ 
a smooth map into the set of all $n\times n$ matrices. 
Let $B(q)$ be the cofactor matrix of $A(q)$. 
We set $A:=A(0)$ and $B:=B(0)$. 
Suppose that 
\begin{equation} \label{b1} 
    \det A = \left.
             \frac{\partial}{\partial q}
             \right|_{q=0} 
              \det A(q) 
           = 0. 
\end{equation}
Then the following lemma holds. 

 \medskip

\begin{Lemma}
Let $X$ be an $n\times n$ matrix 
such that $\tr(XB)\ne 0$. 
Then the following identity holds: 
\begin{eqnarray} \label{b2} 
     \frac{\partial B}{\partial q}(0) 
      & = & \frac{1}{\tr(XB)} 
            \left\{ 
              \tr\left(
               \frac{\partial A}{\partial q}(0) 
                \cdot 
               \left.\frac{\partial Y_t}{\partial t}\right|_{t=0} 
                 \right) 
               \cdot B 
            \right. \\ 
      &   & \qquad 
            \left. 
             -\left.\frac{\partial Y_t}{\partial t}\right|_{t=0} 
               \cdot \frac{\partial A}{\partial q}(0) \cdot B 
             -B\cdot \frac{\partial A}{\partial q}(0) \cdot 
               \left.\frac{\partial Y_t}{\partial t}\right|_{t=0} 
            \right\}, 
            \nonumber
\end{eqnarray}
where $Y_t$ is the cofactor matrix of $A+tX$. 
\end{Lemma}

\begin{Proof}
We set $A_t(q):=A(q)+tX$, 
and denote by $B_t(q)$ its cofactor matrix. 
We have the following Taylor expansions: 
\begin{eqnarray*}
    A_t(q) & = & (A+tX) 
                 +q\frac{\partial A}{\partial q}(0) 
                 +o(q), \\ 
    B_t(q) & = & Y_t 
                 +q\frac{\partial B_t}{\partial q}(0) 
                 +o(q). 
\end{eqnarray*}
Since $A_t(q)B_t(q)=\det A_t(q)\cdot I$, 
we have by taking the first degree terms that 
$$
    \left.
    \frac{\partial}{\partial q}
    \right|_{q=0} 
    \det A_t(q) \cdot I 
     = 
    \frac{\partial A}{\partial q}(0) \cdot Y_t 
     + (A+tX) \cdot \frac{\partial B_t}{\partial q}(0). 
$$
Since 
$$
    \left.
    \frac{\partial}{\partial t}
    \right|_{t=0} 
    \det(A+tX) 
     = \tr(XB) \ne 0, 
$$
$A+tX$ is non-singular around $t=0$. 
Hence we have 
\begin{eqnarray*}
    \frac{\partial B_t}{\partial q}(0) 
     & = & (A+tX)^{-1} 
           \left( 
           \left. 
           \frac{\partial}{\partial q}
           \right|_{q=0} 
           \det A_t(q) \cdot I 
            - 
           \frac{\partial A}{\partial q}(0) \cdot Y_t 
           \right) \\ 
     & = & \frac{\left.
                 \frac{\partial}{\partial q}
                 \right|_{q=0} 
                 \det A_t(q) \cdot Y_t 
                  - 
                 Y_t \cdot 
                 \frac{\partial A}{\partial q}(0) 
                 \cdot Y_t 
                }
                {\det(A+tX)}. 
\end{eqnarray*}
Apply de L'Hospital rule 
to the right-hand side of 
$\frac{\partial B}{\partial q}(0)
=\lim_{t\to 0}\frac{\partial B_t}{\partial q}(0)$. 
Then we get the equality (\ref{b2}). 
\end{Proof}